%% file: arxiv.tex
\documentclass{gaceta}

\usepackage{amsfonts}
\usepackage{amssymb}
\usepackage{amsthm}
\usepackage{amsmath}
\usepackage{epsfig}
\usepackage{psfrag}
\usepackage{url}
\usepackage{mathdots}
\usepackage{xcolor}

\usepackage[spanish]{babel}
\usepackage[T1]{fontenc}
\usepackage[applemac]{inputenc}

\belongstopart{Francisco Santos}

\journame{La Gaceta de la RSME}
\yearofpublication{2010}
\volume{13}
\issuenumber{3}

\def\@begintheorem#1#2{\par\bgroup{\sc #1\ #2. }\it\ignorespaces}
\def\@opargbegintheorem#1#2#3{\par\bgroup{\sc #1\ #2\ (#3). } \it\ignorespaces}
\def\@endtheorem{\egroup}

\newtheorem{teorema}{Teorema}[section]

\newtheorem{lema}[teorema]{Lema}

\newtheorem{definicion}[teorema]{Definición}
\newtheorem{conjetura}[teorema]{Conjetura}

\def\N{\mathbb{N}}

\def\R{\mathbb{R}}

\title{Sobre un contraejemplo a la conjetura de Hirsch}
\author{Francisco Santos\thanks{Financiado parcialmente por el proyecto MTM2008-04699-C03-02 del Ministerio de Ciencia e Innovaci\'on.}} 
\contact{Francisco Santos, Dpto. de Matem\'aticas, Estad\'{\i}stica y Computaci\'on, Universidad de Cantabria, 39005 Santander}
{francisco.santos@unican.es}{http://personales.unican.es/santosf}
\date{}

\begin{document}

\maketitle

\begin{abstract}
Se describe someramente el origen e historia de la Conjetura de Hirsch sobre el di\'ametro posible del grafo de un politopo, as\'{\i} como las ideas principales que llevaron al contraejemplo a la misma recientemente anunciado por el autor.
\end{abstract}

\section{La conjetura y su contexto}
En el a\~{n}o 2000 la revista \emph{Computing in Science and Engineering} pidi\'o a Jack Dongarra y a Francis Sullivan que eligieran los ``10 Algoritmos del Siglo XX''~\cite{Dongarra-Sullivan-topten}; es decir, los algoritmos m\'as influyentes en el desarrollo de la ciencia y la ingenier\'{\i}a del pasado siglo. Uno de los diez elegidos fue el \emph{m\'etodo del s\'{\i}mplice}~\cite{Nash-simplexmethod}.

\subsection*{Programaci\'on lineal y m\'etodo del s\'{\i}mplice\footnote{La mayor\'{\i}a de la literatura sobre 
optimizaci\'on en espa\~nol se refiere a este m\'etodo como el \emph{m\'etodo del s\'{\i}mplex} o el 
\emph{m\'etodo s\'{\i}mplex}, dejando sin traducir la palabra \emph{s\'{\i}mplex} que denota a las celdas de un 
complejo simplicial, o a los politopos con v\'ertices af\'{\i}nmente independientes. Pero puesto que en 
combinatoria geom\'etrica y topol\'ogica es m\'as habitual utilizar  \emph{s\'{\i}mplice}
para referirse a estos objetos, hemos decidido aqu\'{\i}\ aplicar esta palabra tambi\'en al m\'etodo.}}
La programaci\'on lineal es el problema de encontrar el m\'aximo (o m\'{\i}nimo) de una funci\'on lineal en un dominio definido por desigualdades tambi\'en lineales. Naci\'o hacia 1939 con los trabajos del ruso L.~V.~Kantorovitch (1912-1986)~\cite{Kantorovitch}, quien en 1975 recibi\'o por ello el Premio Nobel de Econom\'{\i}a. 
Pero su trabajo no tuvo mucho impacto y qued\'o dormido durante muchos a\~{n}os, por razones estrat\'egicas\footnote{No en vano se trata de la teor\'{\i}a de c\'omo organizar de la mejor manera posible una cantidad limitada de recursos (o defensas) para obtener de ellos el mayor rendimiento (o conseguir los m\'{\i}nimos da\~{n}os).} y pol\'{\i}ticas.\footnote{Aunque la econom\'{\i}a planificada de la URSS podr\'{\i}a parecer el caldo de cultivo perfecto para las ideas de Kantorovitch, los economistas sovi\'eticos eran muy reticentes a que los precios, producci\'on, etc. fueran ``decididos'' por el propio sistema, y en general a cualquier intrusi\'on matem\'atica en su \'ambito.} 
Acabada la guerra, George Dantzig (1914-2005), que por entonces trabajaba en la Oficina de Control Estad\'{\i}stico del ej\'ercito estadounidense, recibi\'o el encargo de desarrollar modelos que permitieran la planificaci\'on (``programaci\'on'', en el lenguaje militar) a gran escala del uso de recursos y horas de trabajo. En 1947 public\'o su m\'etodo del s\'{\i}mplice~\cite{Dantzig-simplex}, que resolv\'{\i}a los problemas de programaci\'on lineal de manera extraordinariamente eficiente, como \'el mismo ilustr\'o con la resoluci\'on de un problema con $70$ variables~\cite{Dantzig-simplex-applied}. 
En la terna de ``padres de la programaci\'on lineal'' suele incluirse tambi\'en al eminente  John von Neumann (1903-1957) quien, tambi\'en en 1947, desarroll\'o la teor\'{\i}a de la dualidad en programas lineales.\footnote{Dantzig relata en~\cite{Dantzig-LPhistory} c\'omo, tras exponerle a von Neumann sus nuevas ideas en una visita privada en Octubre de 1947, \'este le sorprendi\'o con una improvisada ``conferencia de hora y media sobre la teor\'{\i}a matem\'atica de los programas lineales''. 
Ante el asombro de Dantzig, von Neumann confes\'o que lo \'unico que hab\'{\i}a hecho era resumir la teor\'{\i}a de juegos que acababa de desarrollar con Oscar Morgensten, y que su conjetura era que los dos problemas eran equivalentes.}

\begin{figure}
\begin{center}
{\includegraphics[height=4.5 cm]{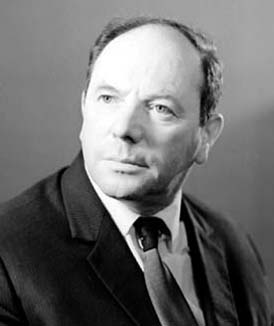}
\quad
\includegraphics[height=4.5 cm]{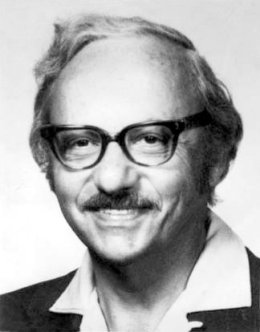}
\quad
\includegraphics[height=4.54 cm]{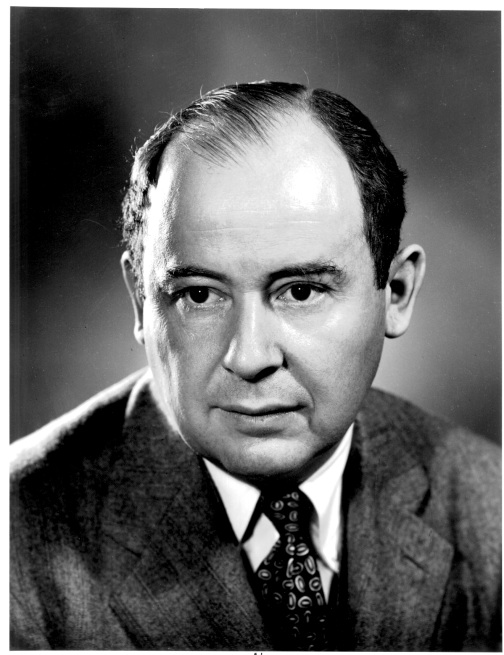}
}
\end{center}
\caption{Kantorovitch, Dantzig y von Neumann, padres de la programaci\'on lineal.}
\end{figure}

La relevancia de la programaci\'on lineal queda expresada, por ejemplo, en las siguientes citas:

\begin{quote}
\emph{Se usa para asignar recursos, planificar producci\'on, organizar turnos de trabajo, dise\~{n}ar carteras de inversi\'on y formular estrategias de mercado, o militares. La versatilidad e impacto econ\'omico de la programaci\'on lineal en el mundo industrial de hoy es verdaderamente admirable.}\,\footnote{\em It is used to allocate resources, plan production, schedule workers, plan investment portfolios and formulate marketing (and military) strategies. The versatility and economic impact of linear programming in today's industrial world is truly awesome.} (Eugene Lawler~\cite{Lawler:sputnik})
\end{quote}

\begin{quote}
\emph{Si hici\'eramos estad\'{\i}sticas sobre qu\'e problema matem\'atico est\'a usando m\'as tiempo de computaci\'on en este momento en el mundo (excluyendo problemas de manejo de bases de datos, como b\'usqueda u ordenaci\'on) la respuesta ser\'{\i}a probablemente la programaci\'on lineal.}\,\footnote{\em If one would take statistics about which mathematical problem is using up most of the computer time in the world, then (not including database handling problems like sorting and searching) the answer would probably be linear programming.} (L\'aszl\'o Lov\'asz~\cite{Lovasz:LP})
\end{quote}

\subsection*{El m\'etodo del s\'{\i}mplice}

 Durante m\'as de 30 a\~{n}os el m\'etodo del s\'{\i}mplice fue el \'unico m\'etodo practicable para resolver grandes problemas de programaci\'on lineal. Sin embargo, a\'un a fecha de hoy no sabemos si es un algoritmo polin\'omico, en el sentido de la teor\'{\i}a de la complejidad. Es decir, no sabemos si un problema de programaci\'on lineal con un n\'umero $d$ de variables y un n\'umero $n$ de restricciones puede ser resuelto mediante el m\'etodo del s\'{\i}mplice en un tiempo que dependa de manera polin\'omica de los par\'ametros $n$ y $d$. La raz\'on fundamental de ese desconocimiento es que no sabemos, dada una regi\'on poli\'edrica de dimensi\'on $d$ y definida por $n$ desigualdades lineales, si su \emph{grafo} tiene \emph{di\'ametro} polin\'omico en los par\'ametros $n$ y $d$. 

Expliquemos esto: Todo dominio (acotado, el caso no acotado es similar) de $\R^d$ definido por restricciones lineales es un \emph{politopo}~\cite{Ziegler:LecturesPolytopes} (lo que en dimensi\'on tres llamamos com\'unmente un poliedro) y tiene v\'ertices, aristas, y caras de diversas dimensiones. Los v\'ertices y aristas forman un grafo, cuyo \emph{di\'ametro} es el m\'aximo n\'umero de aristas que puede ser necesario recorrer para viajar desde un v\'ertice a otro por el camino m\'as corto.
El m\'etodo del s\'{\i}mplice funciona en dos etapas: primero busca un v\'ertice arbitrario del dominio definido por las restricciones y despu\'es va movi\'endose de v\'ertice a v\'ertice, recorriendo en cada paso una arista del politopo y haciendo siempre que el funcional que queremos maximizar aumente. El m\'etodo tiene cierta libertad a la hora de elegir a qu\'e v\'ertice vecino del actual dirigirse y el criterio utilizado para la elecci\'on de uno u otro se llama la ``regla de pivote''. Esta indefinici\'on hace que muchos autores no hablen del m\'etodo del s\'{\i}mplice como un \'unico algoritmo sino como una familia de ellos, cuya complejidad y caracter\'{\i}sticas pueden depender de la regla elegida, que bien puede tener una componente aleatoria.

En todo caso, encontrar el nuevo v\'ertice al que nos queremos mover es computacionalmente sencillo. Pero el n\'umero de veces que hay que hacerlo ser\'a, como m\'{\i}nimo, la distancia del v\'ertice original al v\'ertice \'optimo en el grafo del politopo. Si tenemos mala suerte, esa distancia puede  ser el di\'ametro del grafo.

\subsection*{Las conjeturas de Hirsch y de los $d$ pasos}

Es decir, para poder acotar la complejidad del m\'etodo del s\'{\i}mplice es necesario ser capaces primero de acotar el di\'ametro de los grafos de politopos.\footnote{El rec\'{\i}proco no es cierto; aunque supi\'eramos que el di\'ametro es peque\~{n}o, quedar\'{\i}a el problema de c\'omo hacer que el algoritmo del s\'{\i}mplice encuentre un camino corto.}  Aqu\'{\i} es donde entra la Conjetura de Hirsch:

\begin{conjetura}
[Hirsch 1957]
El di\'ametro (combinatorio) del grafo de un politopo de dimensi\'on $d$ definido por $n$ desigualdades no puede ser nunca mayor que $n-d$.
\end{conjetura}

Esta conjetura fue formulada por Warren M.~Hirsch (1918-2007) en una carta dirigida a Dantzig en 1957. Dantzig la incluy\'o en su libro ``Linear programming and extensions''~\cite{Dantzig-book}, considerado la ``Biblia'' de la programaci\'on lineal. Desde entonces ha atra\'{\i}do la atenci\'on de matem\'aticos tanto puros como aplicados. Sin embargo, m\'as de 50 a\~{n}os despu\'es nuestro conocimiento sobre el problema sigue siendo humillantemente escaso: no se conoce ninguna cota superior polin\'omica para el di\'ametro que se conjetura lineal!! 

Pero algunas cosas s\'{\i}\ que se conocen.\footnote{El lector interesado en los detalles puede consultar el survey cl\'asico de Klee y Kleinschmidt~\cite{Klee-Kleinschmidt} o el que yo he escrito recientemente con Eddie Kim~\cite{Kim-Santos-update}.}
Una de las personas que m\'as contribuy\'o a ello es Victor Klee (1925-2007) profesor de la Universidad de Washington en Seattle, quien, por ejemplo, en 1967 demostr\'o, junto con D.~W.~Walkup, que el estudio de la conjetura en general es equivalente al estudio del caso $n=2d$. (Obs\'ervese que este es, por ejemplo, el caso del $d$-cubo):

\begin{figure}
\begin{center}
{\includegraphics[height=5 cm]{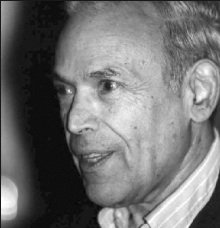}
\qquad
\includegraphics[height=5 cm]{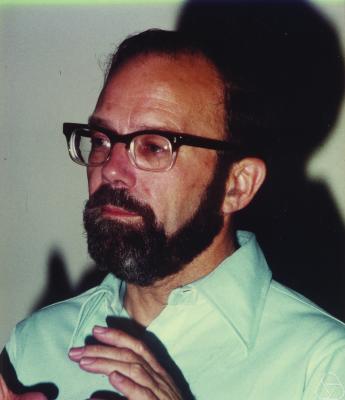}
}
\end{center}
\caption{Hirsch (izquierda), autor de la conjetura que lleva su nombre;  Klee (derecha) demostr\'o, junto con Walkup, su equivalencia con la 
\emph{Conjetura de los $d$ pasos}. Ambos fallecieron en 2007.}
\end{figure}

\begin{teorema}[Teorema de los $d$ pasos~\cite{Klee-Walkup}]
\label{teorema:d-pasos}
Los siguientes enunciados son equivalentes:
\begin{enumerate}
\item El di\'ametro (combinatorio) del grafo de un politopo de dimensi\'on $d$ definido por $n$ desigualdades no puede ser nunca mayor que $n-d$ (Conjetura de Hirsch).
\item El di\'ametro (combinatorio) del grafo de un politopo de dimensi\'on $d$ definido por $2d$ desigualdades no puede ser nunca mayor que $d$ (Conjetura de los $d$ pasos).
\end{enumerate}
\end{teorema}

En el mismo art\'{\i}culo, Klee y Walkup demostraron la Conjetura de los $d$ pasos para $d\le 5$. El caso $d=6$ ha sido verificado computacionalmente hace apenas dos a\~{n}os por Bremner y Schewe~\cite{Bremner:DiameterFewFacets}. De \'el se deduce la Conjetura de Hirsch para  $n\le d+6$, porque la equivalencia a que se refiere el Teorema~\ref{teorema:d-pasos} no es ``para cada $d$\,'' sino, m\'as bien, ``para cada $n-d$\,'': Klee y Walkup demuestran que, si llamamos $H(n,d)$ al di\'ametro m\'aximo que puede tener un politopo de dimensi\'on $d$ con $n$ facetas, se tiene que:
\[
\forall m\in\N, \qquad \max_{d\in \N} \{H(d+m,d)\} = H(2m,m).
\]

\subsection*{Complejidad de la programaci\'on lineal}

En 1979 y 1984, L.~G.~Khachiyan~\cite{Khachiyan} y N.~Karmarkar~\cite{Karmarkar} encontraron sendos algoritmos polin\'omicos para resolver la programaci\'on lineal: el \emph{m\'etodo del elipsoide} y el \emph{m\'etodo de puntos interiores}. El primero, aunque supuso una revoluci\'on desde el punto de vista te\'orico~\cite{Lawler:sputnik,Lovasz:LP}, nunca tuvo muchas implicaciones pr\'acticas por su dificultad de implementaci\'on y sus pobres resultados. El segundo s\'{\i}, pero a\'un hoy el m\'etodo del s\'{\i}mplice compite con \'el en utilizaci\'on, y suscita inter\'es (y perplejidad) por varias razones. Entre otras:

\begin{itemize}

\item Aunque no sabemos si el m\'etodo del s\'{\i}mplice es polin\'omico, en la pr\'actica es un m\'etodo muy r\'apido. Por ejemplo, en un survey sobre programaci\'on lineal Todd~\cite{Todd:LP} dice: \emph{``El m\'etodo del s\'{\i}mplice, t\'{\i}picamente, encuentra el \'optimo en como mucho $2m$ \'o $3m$ pasos'',}\,\footnote{\em ``The (primal) simplex method typically requires at most $2m$ to $3m$ pivots to attain optimality''.
} donde $m=n-d$ representa el n\'umero de desigualdades menos la dimensi\'on.\footnote{En su \emph{forma est\'andar}, el dominio de un problema de programaci\'on lineal se presenta mediante $n$ variables no negativas restringidas por $m$ ecuaciones lineales. La no-negatividad de cada variable produce una faceta, y las ecuaciones bajan la dimensi\'on del dominio a $n-m$.}
En una recensi\'on que acaba de ser publicada en el \emph{Bulletin of the American Math. Society}~\cite{Todd:Dantzig}, Todd repite esa estimaci\'on y a\~{n}ade: \emph{``El m\'etodo del s\'{\i}mplice sigue siendo, si no el m\'etodo a elegir, s\'{\i}\ uno de los m\'etodos a elegir, a la misma altura de, y en ciertas clases de problemas superior a, aproximaciones m\'as modernas''.}\,\footnote{
\em ``The simplex method has remained, if not the method 
of choice, a method of choice, usually competitive with, and on some classes of 
problems superior to, the more modern approaches''.}

Esto se ha considerado tradicionalmente evidencia a favor de la Conjetura de Hirsch o, al menos, a favor de una versi\'on ``polin\'omica'' de la misma. Pero tambi\'en puede ser debido a otras razones. Por ello, hay una vasta literatura dedicada a estudiar qu\'e pasa en el m\'etodo del s\'{\i}mplice para politopos aleatorios, o para reglas de pivote aleatorias en politopos fijos, o cuando se permite perturbar ligeramente los coeficientes de las desigualdades que los definen, etc.~\cite{Borgwardt,Kalai:A-subexponential-randomized,Klee-Minty,Matousek:A-subexponential-bound,Megiddo:LPinLT,Megiddo:SurveyLPCompelxity,Smale-AverageLP,Spielman:WhySimplexUsually,Vershynin:BeyondHirsch}. 

\item Los m\'etodos de Khachiyan y Karmarkar son polin\'omicos cuando la complejidad se mide en el \emph{modelo de bits}. O sea, el ``tama\~{n}o del input'' no son los par\'ametros $n$ y $d$ sino el tama\~{n}o en bits de la matriz que define a nuestro politopo.\footnote{B\'asicamente, eso da $nd \log(K)$ donde $K$ es la cota m\'axima para los coeficientes de la matriz.} Aunque el modelo de bits es muy bueno para el an\'alisis te\'orico, se considera que el \emph{modelo real}, introducido por Blum, Shub y Smale hacia 1990~\cite{BCSS,BSS}, describe mejor la eficiencia pr\'actica de un algoritmo. La pregunta de si existe un algoritmo para programaci\'on lineal que sea polin\'omico en el modelo real est\'a abierta, y fue incluida en el a\~{n}o 2000 por Steven Smale en su lista de \emph{``Problemas matem\'aticos para el siglo que viene''}~\cite{Smale}. Si se encontrase una regla de pivote polin\'omica para el m\'etodo del s\'{\i}mplice, la respuesta a la pregunta de Smale ser\'{\i}a autom\'aticamente afirmativa.

\end{itemize}

\section{El contraejemplo}

El 10 de Mayo envi\'e lo siguiente como resumen para mi conferencia en el congreso \emph{The Mathematics of Klee and Gr\"unbaum: 100 years in Seattle}~\cite{100-years-in-Seattle}:

\begin{quote}
\emph{S\'olo he estado en Seattle una vez, en enero de 2002, para impartir un coloquio en la U.~de Washington. Aunque Victor Klee ya estaba jubilado --ten\'{\i}a 76 a\~{n}os-- vino al Departamento de Matem\'aticas para charlar conmigo. Tuvimos una amena conversaci\'on en el transcurso de la cual me pregunt\'o: ¿Por qu\'e no intentas refutar la Conjetura de Hirsch? [\dots] Esta charla da respuesta a esa pregunta. En ella describir\'e la construcci\'on de un politopo de dimensi\'on $43$ con $86$ facetas y di\'ametro mayor que $43$. La prueba se basa en una generalizaci\'on del Teorema de los $d$ pasos de Klee y Walkup.}\,%
\footnote{
\emph
{I have been in Seattle only once, in January 2002, when I visited to give a colloquium talk at UW. Although Victor Klee was already retired---he was 76 years old---he came to the Department of Mathematics to talk to me. We had a nice conversation during which he asked: Why don't you try to disprove the Hirsch Conjecture? [\dots] This talk is the answer to that question. I will describe the construction of a $43$-dimensional polytope with $86$ facets and diameter bigger than $43$. The proof is based on a generalization of the $d$-step Theorem of Klee and Walkup.}%
}
\end{quote}

Ese mismo d\'{\i}a Gil Kalai public\'o la noticia en su muy visitado blog~\cite{kalaiblog} y la entrada de ``Hirsch conjecture'' en la Wikipedia fue actualizada para hacerse eco de ella. Incluso, alguien escribi\'o un soneto al respecto~\cite{readwriteweb}.

\medskip
El contraejemplo que anuncio en el resumen tiene dos ingredientes: una generalizaci\'on del Teorema de los $d$ pasos, y la construcci\'on expl\'{\i}cita de cierto politopo de dimensi\'on cinco con $48$ facetas y unas propiedades concretas.

\subsection*{El Teorema generalizado de los $d$ pasos}

El \emph{Teorema generalizado de los $d$ pasos} se refiere a una clase muy particular de politopos, que llamaremos \emph{husos}:

\begin{definicion}
Un \emph{huso} es un politopo que posee dos v\'ertices con la propiedad de que toda \emph{faceta} (cara de codimensi\'on uno) contiene a alguno de ellos. A esos dos v\'ertices los llamamos los \emph{extremos} del huso.
\end{definicion}

\begin{figure}[htb]
\begin{center}
\resizebox{0.6\linewidth}{!}{
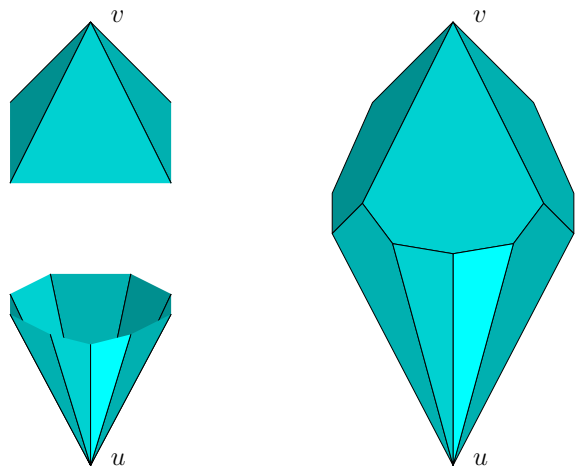
}
\caption{Un huso se obtiene al intersecar dos conos no acotados. 
}
\label{fig:huso}
\end{center}
\end{figure}

Dicho de otro modo, un huso es la intersecci\'on de dos conos poli\'edricos no acotados, cuyos \'apices son los extremos del huso, como en la Figura~\ref{fig:huso}.

Obs\'ervese que las facetas de un politopo se corresponden exactamente con las desigualdades que lo definen (cada faceta est\'a contenida en el hiperplano donde una de las desigualdades se satisface con igualdad). Por tanto, el n\'umero de facetas  es el n\'umero de desigualdades en el problema de la programaci\'on lineal (salvo que haya desigualdades superfluas). 

\begin{teorema}
[Teorema generalizado de los $d$ pasos]
\label{teorema:d-step-gen}
Los siguientes enunciados son equivalentes:
\begin{enumerate}
\item El di\'ametro (combinatorio) del grafo de un politopo de dimensi\'on $d$ definido por $n$ desigualdades no puede ser nunca mayor que $n-d$ (Conjetura de Hirsch).
\item En todo huso de dimensi\'on $d$ se puede ir de un extremo al otro atravesando como mucho $d$ aristas.
\end{enumerate}
\end{teorema}

Lo importante aqu\'{\i}\ es que
el n\'umero de facetas del huso ya no aparece en el enunciado, como ocurr\'{\i}a en el Teorema de los $d$ pasos original.

La demostraci\'on del Teorema~\ref{teorema:d-step-gen} se basa en el siguiente lema:

\begin{lema}[Lema del Huso]
\label{lema:huso}
Si $P$ es un huso de dimensi\'on $d$, con $m$ facetas, y en el que hacen falta $k$ pasos para ir de un extremo a otro, y si $m>2d$, entonces existe otro huso de dimensi\'on $d+1$, con $m+1$ facetas, y en el que hacen falta $k+1$ pasos para ir de un extremo a otro.
\end{lema}

El Lema del Huso nos da la implicaci\'on de (2) a (1) en el Teorema~\ref{teorema:d-step-gen} por inducci\'on sobre $m-2d$, donde $m$ es el n\'umero de facetas de nuestro huso: Si $m\le 2d$ entonces no hay nada que demostrar (el huso es ya un contraejemplo a la Conjetura de Hirsch) y si $m > 2d$, aplicando $m-2d$ veces el lema obtenemos un politopo (un huso, aunque esto ya no es importante) de dimensi\'on $d + (m-2d) = m-d$, con n\'umero de facetas $2m-2d$, y en el que hacen falta $k + (m-2d) > d + (m-2d) = m-d$ pasos para ir de un extremo a otro. O sea, hemos violado la Conjetura de los $d$ pasos ``tradicional'' y, en particular, la Conjetura de Hirsch.

Es f\'acil demostrar que en todo huso $3$-dimensional bastan tres aristas para ir de un extremo $u$ al otro $v$: los v\'ertices y aristas no incidentes ni a $u$ ni a $v$ forman un ciclo que necesariamente contiene alg\'un v\'ertice adyacente a $u$, alguno adyacente a $v$, y alguna arista que conecta v\'ertices de los dos tipos (ver Figura~\ref{fig:huso}).
%
Pero en dimensi\'on superior no ocurre as\'{\i}. El contraejemplo a la Conjetura de Hirsch se obtiene a partir del siguiente huso:

\begin{teorema}
\label{teorema:huso}
Existe un huso de dimensi\'on $5$ con $48$ facetas y en el que hacen falta $6$ pasos para ir de un extremo al otro.
\end{teorema}

Como hemos de aplicar el Lema del Huso $n-2d = 38$ veces, el resultado final es un politopo de dimensi\'on $43$ y con $86$ facetas. 

\subsection*{Husos y prismatoides}

El Teorema~\ref{teorema:huso} ha sido verificado computacionalmente por Julian Pfeifle y Eddie Kim con la ayuda del software \url{polymake}~\cite{polymake}\footnote{Pueden verse algunos detalles en \url{http://personales.unican.es/santosf}}. Sin embargo, el art\'{\i}culo~\cite{Santos:Hirsch-counter} contiene una demostraci\'on del mismo, previa a la verificaci\'on computacional y que, de hecho, explica c\'omo se lleg\'o a la construcci\'on del huso en cuesti\'on. Inclu\'{\i}mos aqu\'{\i} algunas de las ideas, que permiten entender un huso de dimensi\'on cinco como un objeto puramente $3$-dimensional.

Quiz\'a por ``deformaci\'on profesional'', puesto que he trabajado mucho con triangulaciones y complejos simpliciales~\cite{triang-book}, mi tendencia es siempre a dualizar la Conjetura de Hirsch: no pienso en c\'omo ir de un v\'ertice a otro en un politopo a trav\'es de sus aristas, sino en c\'omo ir de una faceta a otra en pasos que consisten en cruzar caras de codimensi\'on dos (\emph{crestas}). Esto hace el problema mucho m\'as combinatorio. Lo primero que debemos entender es en qu\'e se convierten los husos al dualizarlos:

\begin{definicion}
Un \emph{prismatoide} es un politopo cuyos v\'ertices est\'an todos contenidos en la uni\'on de dos facetas, a las que llamaremos \emph{bases}.
\end{definicion}

El nombre proviene, por supuesto, de que un prisma es un caso particular de prismatoide. Animamos al lector a que vuelva a pensar qu\'e significa la propiedad de los $d$ pasos para un prismatoide de dimensi\'on tres: al principio nos encontramos en una de las dos bases, digamos la de abajo; en un primer paso vamos a una de las facetas que comparten arista con dicha base. En un segundo paso nos trasladamos a una que comparte arista con la base de arriba, y en un tercero vamos a la base de arriba.

\subsection*{Prismatoides y mapas geod\'esicos superpuestos}


Sea $Q$ un prismatoide de dimensi\'on $d$, con bases $Q^+$ y $Q^-$. Las bases son a su vez politopos de dimensi\'on $d-1$, y no hay p\'erdida de generalidad en suponerlas paralelas (eso se puede conseguir mediante una transformaci\'on proyectiva, que no altera la combinatoria de $Q$). Sea $H$ un hiperplano paralelo a $Q^+$ y $Q^-$ e intermedio entre ambos. No es dif\'{\i}cil demostrar entonces que:

\begin{itemize}
\item La secci\'on $Q\cap H$ es suficiente para recuperar toda la combinatoria de $Q$. (Ver Figura~\ref{fig:prismatoide}.)
\begin{figure}
\begin{center}
\resizebox{0.97\linewidth}{!}{
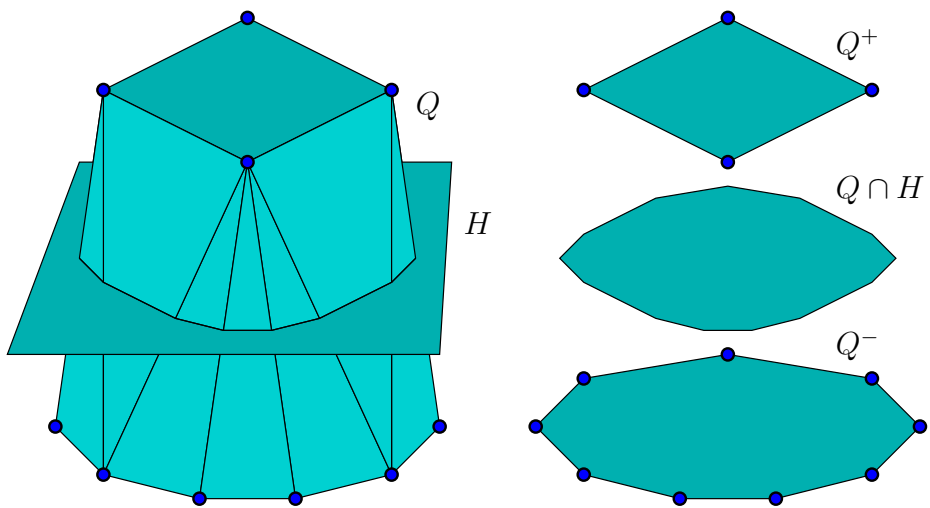
}
\end{center}
\caption{Dos $d-1$-politopos $Q^+$ y $Q^+-$, el $d$-prismatoide $Q$ que forman, y su intersecci\'on con un hiperplano intermedio}
\label{fig:prismatoide}
\end{figure}

\item Combinatoriamente, $Q \cap H$ es la suma de Minkowski de $Q^+$ y $Q^-$. (Ver Figura~\ref{fig:suma}.)
\begin{figure}
\begin{center}
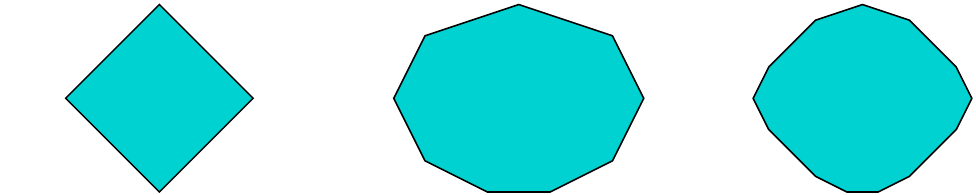
\end{center}
\caption{Toda secci\'on paralela del prismatoide es la suma de Minkowski (ponderada) de las bases $Q^+$ y $Q^-$}
\label{fig:suma}
\end{figure}
%
En particular, el abanico normal de $Q\cap H$ es la superposici\'on (``refinamiento com\'un'') de los abanicos normales de $Q^+$ y $Q^-$.
\end{itemize}

Recordemos lo que es el \emph{abanico normal} de un politopo $P$: en cada punto $x$ del borde del politopo se puede definir el cono normal exterior, que no es m\'as que el conjunto de funcionales lineales cuyo m\'aximo en $P$ se alcanza en $x$. Puntos en (el interior relativo de) una misma cara de $P$ tienen el mismo cono normal, as\'{\i}\ que $P$ nos da una descomposici\'on del espacio vectorial $(\R^d)^*$ de todos los funcionales lineales en una familia finita de conos con \'apice en el origen y que, adem\'as, forman un complejo en el sentido habitual (sus intersecciones son siempre en caras completas). Este complejo es el abanico normal de $P$.

Ahora bien, los abanicos normales de los politopos $Q^+$ y $Q^-$, que est\'an formados por conos lineales de dimensi\'on $d-1$, podemos intersecarlos con la esfera unidad $S^{d-2}$ sin perder informaci\'on sobre ellos. Por tanto, construir un prismatoide de dimensi\'on cinco es b\'asicamente lo mismo que construir dos \emph{mapas geod\'esicos} (descomposiciones celulares en las que todas las celdas son intersecci\'on de hemisferios) en la 3-esfera $S^3$. El siguiente lema nos dice qu\'e han de cumplir esos dos mapas para que el prismatoide nos d\'e un contraejemplo a la Conjetura de Hirsch:

\begin{lema}
 Sea $Q$ un prismatoide de dimensi\'on d con bases paralelas $Q^+$ y $Q^-$ y sean $G^+$ y $G^-$ los mapas geod\'esicos que se obtienen al intersecar con $S^{d-2}$ los abanicos normales de $Q^+$ y $Q^-$. Entonces, las siguientes propiedades son equivalentes:
\begin{enumerate}
\item Se puede ir de $Q^+$ a $Q^-$ en Q en a lo m\'as $d$ pasos.
\item Al superponer $G^+$ y $G^-$, se puede ir de un v\'ertice de $G^+$ a uno de $G^-$ atravesando a lo m\'as $d-2$ aristas del refinamiento com\'un as\'{\i}\ obtenido.
\end{enumerate}
\end{lema}

\begin{figure}[htb]
\begin{center}
\includegraphics[scale=0.8]{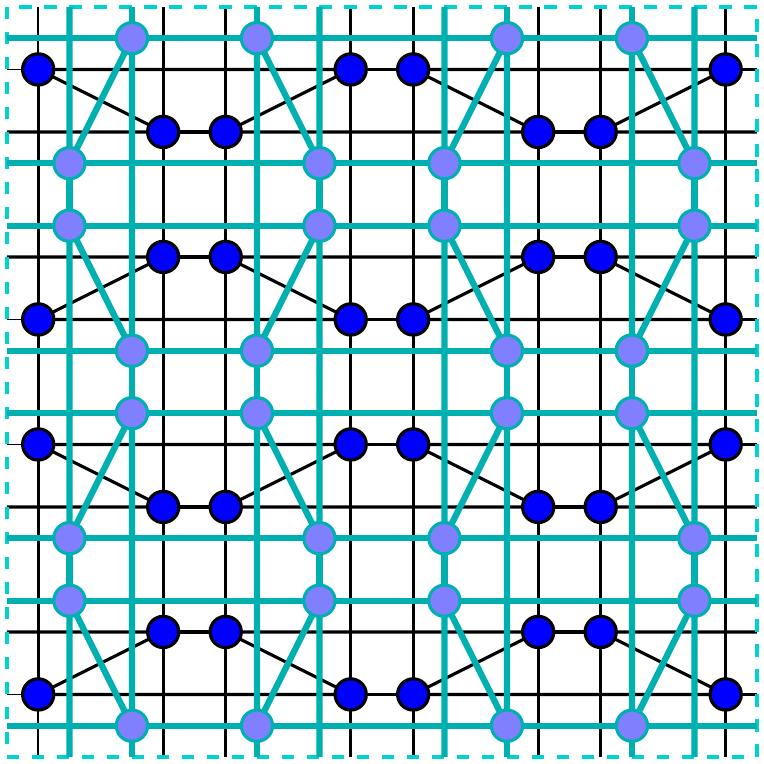}
\caption{Los mapas geod\'esicos de las bases $Q^+$ (oscuro) y $Q^-$ (claro) del prismatoide dual al huso del Teorema~\ref{teorema:huso}. Los mapas est\'an dibujados en un mismo toro, pero han de entenderse en dos toros paralelos en la $3$-esfera}
\label{fig:Q}
\end{center}
\end{figure}

Es decir, construir un prismatoide en dimensi\'on cinco que viole la propiedad de los $d$ pasos es lo mismo que construir una pareja de mapas geod\'esicos en la esfera $S^3$ que al supeponerlos violen la ``propiedad de los $d-2$ pasos'' expresada en este \'ultimo lema. Esto es lo que hace mi construcci\'on, usando en cierto modo para ello la bien conocida \emph{fibraci\'on de Hopf} de la $3$-esfera: los dos mapas tienen cada uno su ``parte interesante'' en un toro, y esos dos toros se colocan uno paralelo al otro, y de manera adecuada, en la esfera. Los mapas concretos utilizados, dibujados en un toro, se muestran en la Figura~\ref{fig:Q}. Pero los detalles los dejamos para que el lector interesado los consulte en el art\'{\i}culo~\cite{Santos:Hirsch-counter}.

\section{¿Y ahora qu\'e?}

En un survey sobre la Conjetura de Hirsch~\cite{Klee-Kleinschmidt} publicado en 1987, Klee y Kleinschmidt dicen:
\emph{Encontrar un contraejemplo a la conjetura no ser\'a m\'as que un peque\~{n}o primer paso en la l\'{\i}nea de investigaci\'on de la misma.}\,\footnote{
\em Finding a counterexample will be merely a small first step in the line of investigation related to the conjecture.}
Aunque para dar ese ``peque\~{n}o paso'' hayan sido necesarios $53$ a\~{n}os desde el enunciado de la conjetura y $23$ desde que se escribieron esas palabras,
suscribo totalmente la frase. Mi contraejemplo se puede convertir, mediante construcciones cl\'asicas,  en una familia inifinita de contraejemplos a la Conjetura de Hirsch en la que el di\'ametro de los politopos construidos crece, b\'asicamente, como $1.02 n$. Es decir, se viola la conjetura, que era $n-d$, pero el di\'ametro obtenido no deja de ser lineal y no nos dice mucho sobre el problema de fondo. 

Quiz\'a, por tanto, m\'as significativo que el contraejemplo en s\'{\i}\ es el \emph{Teorema generalizado de los $d$ pasos}, que abre una nueva l\'{\i}nea de ataque al problema. En este sentido, ser\'{\i}a interesante estudiar los siguientes problemas:

\begin{enumerate}
\item ¿Existen husos/prismatoides para $d=4$ en los que no sea posible ir de un extremo/base el otro en $d$ pasos? Como ya hemos dicho, en dimensi\'on tres es f\'acil demostrar que no existen y en dimensi\'on cinco se tiene el del Teorema~\ref{teorema:huso}.

\item Mi prismatoide, con $24$ v\'ertices en cada base, con toda seguridad no es el m\'as peque\~{n}o posible. Su n\'umero de v\'ertices proviene de que la construcci\'on es altamente sim\'etrica, como se desprende de la Figura~\ref{fig:Q}. ¿Existen $5$-prismatoides con, digamos, $20$ v\'ertices y que necesiten m\'as de cinco pasos para ir de una base a la otra? Si la respuesta es s\'{\i} tendr\'{\i}amos contraejemplos a la Conjetura de Hirsch en dimensi\'on $15$ en vez de $43$.

\item Un desaf\'{\i}o computacional: ¿c\'omo construir y verificar expl\'{\i}citamente el politopo de dimensi\'on 43 que viola la Conjetura de Hirsch? Nuestra demostraci\'on se basa en aplicar $38$ veces el Lema del Huso. Pero cada aplicaci\'on del lema duplica (m\'as o menos) el n\'umero de v\'ertices del huso (o de facetas del prismatoide dual), con lo cual al final tendremos un politopo de dimensi\'on $43$ con m\'as o menos $2^{40}$ (un bill\'on) de v\'ertices.\footnote{Si a alg\'un lector le parece imposible que un politopo definido por una matriz de ``s\'olo'' $86$ desigualdades y $43$ variables tenga ese n\'umero de v\'ertices, que piense en el cubo de dimensi\'on 43, que tiene exactamente $2^{43}$ v\'ertices. Es en ejemplos como \'este donde se aprecia la importancia de la distinci\'on entre ``polin\'omico'' y ``exponencial''.} 

Para empeorar las cosas, la demostraci\'on del Lema del Huso necesita de una ``perturbaci\'on'' de los coeficientes. Por tanto,
aunque los coeficientes del huso/prismatoide original son extremadamente sencillos (n\'umeros enteros entre $-5$ y $5$) esa torre de $38$ perturbaciones da lugar o bien a coeficientes racionales de un tama\~{n}o inmenso si se decide usar aritm\'etica exacta, o bien a problemas muy sofisticados de c\'alculo num\'erico si se decide usar aproximaciones num\'ericas.
\end{enumerate}

\vskip.5cm
\noindent {Francisco Santos}\newline
\emph{Departamento de Matem\'aticas, Estad\'istica y Computaci\'on}\newline
\emph{Universidad de Cantabria, E-39005 Santander, Spain}\newline
\emph{email: }\url{francisco.santos@unican.es}\newline
\emph{web: }\url{http://personales.unican.es/santosf/}

\end{document}

%% file: spindle.pdf_t
\begin{picture}(0,0)%
\includegraphics{spindle.pdf}%
\end{picture}%
\setlength{\unitlength}{3947sp}%
\begingroup\makeatletter\ifx\SetFigFont\undefined%
\gdef\SetFigFont#1#2#3#4#5{%
  \reset@font\fontsize{#1}{#2pt}%
  \fontfamily{#3}\fontseries{#4}\fontshape{#5}%
  \selectfont}%
\fi\endgroup%
\begin{picture}(4224,3528)(289,-3130)
\put(1051,-3061){\makebox(0,0)[lb]{\smash{{\SetFigFont{12}{14.4}{\rmdefault}{\mddefault}{\updefault}{$u$}%
}}}}
\put(3751,-3061){\makebox(0,0)[lb]{\smash{{\SetFigFont{12}{14.4}{\rmdefault}{\mddefault}{\updefault}{$u$}%
}}}}
\put(3751,239){\makebox(0,0)[lb]{\smash{{\SetFigFont{12}{14.4}{\rmdefault}{\mddefault}{\updefault}{$v$}%
}}}}
\put(1051,239){\makebox(0,0)[lb]{\smash{{\SetFigFont{12}{14.4}{\rmdefault}{\mddefault}{\updefault}{$v$}%
}}}}
\end{picture}%

%% file: prismatoid.pdf_t
\begin{picture}(0,0)%
\includegraphics{prismatoid.pdf}%
\end{picture}%
\setlength{\unitlength}{3947sp}%
\begingroup\makeatletter\ifx\SetFigFont\undefined%
\gdef\SetFigFont#1#2#3#4#5{%
  \reset@font\fontsize{#1}{#2pt}%
  \fontfamily{#3}\fontseries{#4}\fontshape{#5}%
  \selectfont}%
\fi\endgroup%
\begin{picture}(5765,3104)(-4811,-2963)
\put(376,-136){\makebox(0,0)[lb]{\smash{{\SetFigFont{12}{14.4}{\rmdefault}{\mddefault}{\updefault}{$Q^+$}%
}}}}
\put(376,-2011){\makebox(0,0)[lb]{\smash{{\SetFigFont{12}{14.4}{\rmdefault}{\mddefault}{\updefault}{$Q^-$}%
}}}}
\put(376,-1036){\makebox(0,0)[lb]{\smash{{\SetFigFont{12}{14.4}{\rmdefault}{\mddefault}{\updefault}{$Q\cap H$}%
}}}}
\put(-1949,-1261){\makebox(0,0)[lb]{\smash{{\SetFigFont{12}{14.4}{\rmdefault}{\mddefault}{\updefault}{$H$}%
}}}}
\put(-2249,-511){\makebox(0,0)[lb]{\smash{{\SetFigFont{12}{14.4}{\rmdefault}{\mddefault}{\updefault}{$Q$}%
}}}}
\end{picture}%

%% file: sum.pdf_t
\begin{picture}(0,0)%
\includegraphics{sum.pdf}%
\end{picture}%
\setlength{\unitlength}{3947sp}%
\begingroup\makeatletter\ifx\SetFigFont\undefined%
\gdef\SetFigFont#1#2#3#4#5{%
  \reset@font\fontsize{#1}{#2pt}%
  \fontfamily{#3}\fontseries{#4}\fontshape{#5}%
  \selectfont}%
\fi\endgroup%
\begin{picture}(4677,924)(3736,-298)
\put(5101, 89){\makebox(0,0)[lb]{\smash{{\SetFigFont{14}{16.8}{\rmdefault}{\mddefault}{\updefault}{$+\ \frac{1}{2}$}%
}}}}
\put(3751, 89){\makebox(0,0)[lb]{\smash{{\SetFigFont{14}{16.8}{\rmdefault}{\mddefault}{\updefault}{$\frac{1}{2}$}%
}}}}
\put(6976, 89){\makebox(0,0)[lb]{\smash{{\SetFigFont{14}{16.8}{\rmdefault}{\mddefault}{\updefault}{$=$}%
}}}}
\end{picture}%

%% file: arxiv.bbl
\begin{thebibliography}{}


\bibitem{100-years-in-Seattle} ``The Mathematics of Klee and Gr\"unbaum: 100 years in Seattle'', 28--30 de  Julio de 2010, University of Washington, Seattle. Comit\'e organizador Fred Holt, Isabella Novik, Rekha Thomas and Gordon Williams. \url{https://sites.google.com/a/alaska.edu/kleegrunbaum/}







\bibitem{BCSS}
L.~Blum, F.~Cucker, M.~Shub, and S.~Smale,
Complexity and real computation, Springer-Verlag, 1997. 

\bibitem{BSS}
L.~Blum, M.~Shub, and S.~Smale, On a Theory of Computation and Complexity over the Real Numbers: $NP$-completeness, Recursive Functions and Universal Machines. \emph{Bull. Amer. Math. Soc.} {\bf 21} (1), (1989), 1--46.

\bibitem{Borgwardt}
K.~H.~Borgwardt, 
The Average Number of Steps Required by the 
Simplex Method Is Polynomial.
\emph{Zeitschrift fur Operations Research}, 26:157--177, 1982.



\bibitem{Bremner:DiameterFewFacets}
D.~Bremner and L.~Schewe.
Edge-graph diameter bounds for convex polytopes with few facets.
preprint \url{arXiv:0809.0915}, 9 pages, September 2008. Accepted in~\emph{Experimental Mathematics}.



\bibitem{Dantzig-simplex}
G.~B.~Dantzig,  Maximization of a linear function of variables subject to linear inequalities. \emph{Activity Analysis of Production and Allocation}, pp. 339--347. Cowles Commission Monograph No. 13. John Wiley \& Sons, Inc., New York, N. Y.; Chapman \& Hall, Ltd., London, 1951.

\bibitem{Dantzig-simplex-applied}
G.~B.~Dantzig,  Application of the simplex method to a transportation problem. Activity Analysis of Production and Allocation, pp. 359--373. Cowles Commission Monograph No. 13. John Wiley \& Sons, Inc., New York, N. Y.; Chapman \& Hall, Ltd., London, 1951.

\bibitem{Dantzig-book}
G.~B.~Dantzig,
\newblock Linear programming and extensions,
\newblock {Princeton University Press}, {1963}.
Reprinted in the series~\emph{Princeton Landmarks in Mathematics}, Princeton University Press, 1998.


\bibitem{Dantzig-LPhistory}
G.~B.~Dantzig,
\newblock Linear programming,
\newblock in \emph{History of Mathematical Programming: A Collection of Personal Reminiscences}, 
J. K. Lenstra, A. H. G. Rinnooy Kan, and A. Schrijver (eds.), Elsevier Science Publishers B.V., 1991, pp.19--31.





\bibitem{triang-book}
J.~A.~De Loera, J.~Rambau, F.~Santos, 
\emph{Triangulations: Structures for Algorithms and Applications},
Springer-Verlag, 2010. ISBN: 978-3-642-12970-4.





\bibitem{Dongarra-Sullivan-topten}
J.~Dongarra and F.~Sullivan.
Guest Editors' Introduction: The Top 10 Algorithms.
\emph{Comput. Sci. Eng.} 2, 22 (2000), 2 pages.




%


\bibitem{polymake} 
E.~Gawrilow, M.~Joswig.
Polymake: A software package for analyzing convex polytopes. 
In \emph{Polytopes---combinatorics and computation (Oberwolfach, 1997)},
DMV Sem.~\textbf{29}, Birkh\"auser, 2000, pp. 43--73.





\bibitem{readwriteweb}
Curt Hopkins.
The Hirsch Conjecture, Disproved.
\emph{ReadWriteWeb}, \url{http://www.readwriteweb.com/archives/the_hirsch_conjecture_disproved.php}



\bibitem{Kalai:A-subexponential-randomized}
G.~Kalai.
\newblock A subexponential randomized simplex algorithm.
\newblock In {\em Proceedings of the 24th annual ACM symposium on the Theory of
  Computing}, pages 475--482, Victoria, 1992. ACM Press.

\bibitem{kalaiblog} G. Kalai. Online blog~\url{http://gilkalai.wordpress.com}. 


\bibitem{Kantorovitch}
L.~V.~Kantorovitch.
A new method of solving of some classes of extremal problems. 
\emph{C. R. (Doklady) Acad. Sci. URSS (N.S.)} 28, (1940). 211--214.

\bibitem{Karmarkar}
N.~Karmarkar.
\newblock A new polynomial time algorithm for linear programming.
\newblock {\em Combinatorica}, 4(4):373--395, 1984.

\bibitem{Khachiyan}
L.~G.~Ha\v{c}ijan. 
\newblock A polynomial algorithm in linear programming. (in Russian)
\newblock \emph{Dokl. Akad. Nauk SSSR}, 244(5):1093--1096, 1979.
 

\bibitem{Kim-Santos-update}
E.~D. Kim, and F.~Santos.
\newblock An update on the Hirsch conjecture,
{\em Jahresbericht der Deutschen Mathematiker-Vereinigung (DMV)}, to appear 2010.
DOI: 10.1365/s13291-010-0001-8


\bibitem{Klee-Kleinschmidt}
V.~Klee, P.~Kleinschmidt, 
The $d$-Step Conjecture and Its Relatives,
\emph{Mathematics of Operations Research}, 12(4):718--755, 1987. 

\bibitem{Klee-Minty}
V.~Klee, G.~J.~Minty, How good is the simplex algorithm?, in \emph{Inequalities, III (Proc. Third Sympos., Univ. California, Los Angeles, Calif., 1969; dedicated to the memory of Theodore S. Motzkin)},  Academic Press, New York, 1972, pp. 159--175.

\bibitem{Klee-Walkup}
V.~Klee and D.~W. Walkup.
\newblock The $d$-step conjecture for polyhedra of dimension $d < 6$.
\newblock {\em Acta Math.}, 133:53--78, 1967.



\bibitem{Lawler:sputnik}
E.~L.~Lawler.
The great mathematical sputnik of 1979.
\emph{The Mathematical Intelligencer}
Volume 2, Number 4 (1980), 191--198.

\bibitem{Lovasz:LP}
L.~Lov\'asz.
A new linear programming algorithm---better or worse than the simplex method? 
\emph{Math. Intelligencer} 2, no. 3 (1979/80), 141--146.


\bibitem{Matousek:A-subexponential-bound}
J.~Matou\v{s}ek, M.~Sharir, and E.~Welzl.
\newblock A subexponential bound for linear programming.
\newblock In {\em Proceedings of the 8th annual symposium on Computational
  Geometry}, pages 1--8, 1992.

\bibitem{Megiddo:LPinLT}
N.~Megiddo. 
Linear programming in linear time when the dimension is fixed. 
\emph{J. Assoc. Comput. Mach.}, 31(1):114--127, 1984.

\bibitem{Megiddo:SurveyLPCompelxity}
N.~Megiddo. 
On the complexity of linear programming. 
In: \emph{Advances in economic theory: Fifth world congress}, T.~Bewley, ed.
Cambridge University Press, Cambridge, 1987, 225-268.


\bibitem{Nash-simplexmethod}
J.~C.~Nash.
The (Dantzig) Simplex Method for Linear Programming
\emph{Comput. Sci. Eng.} 2, 29 (2000), 3 pages.




\bibitem{Santos:Hirsch-counter}
F.~Santos.
A counter-example to the Hirsch Conjecture.
Preprint~\url{arXiv:1006.2814}, Junio de 2010.

\bibitem{Smale-AverageLP}
S.~Smale, 
On the Average Number of Steps of the Simplex Method  of Linear Programming. 
\emph{Mathematical Programming}, 27: 241-62, 1983.

\bibitem{Smale} S. Smale, Mathematical problems for the next century. {\em Mathematics: frontiers and perspectives}, pp. 271--294, American Mathematics Society, Providence, RI (2000).

\bibitem{Spielman:WhySimplexUsually}
D.~A. Spielman and S.~Teng.
\newblock Smoothed analysis of algorithms: Why the simplex algorithm usually
  takes polynomial time.
\newblock {\em J. ACM}, 51(3):385--463, 2004.



\bibitem{Todd:LP}
M.~J.~Todd.
The many facets of linear programming 
\emph{Math. Program., Ser. B} 91:3, (2001).

\bibitem{Todd:Dantzig}
M.~J.~Todd.
Book review on ``The basic George B. Dantzig'' (by Richard W. Cottle, Stanford University Press, Stanford, California, 2003). 
\emph{Bull. Amer. Math. Soc. (N. S.)}, en prensa. Publicado online el 19 de Mayo de 2010, DOI: 10.1090/S0273-0979-2010-01303-3.


\bibitem{Vershynin:BeyondHirsch}
R.~Vershynin.
\newblock Beyond {H}irsch conjecture: walks on random polytopes and smoothed
  complexity of the simplex method.
\newblock In {\em IEEE Symposium on Foundations of Computer Science},
  volume~47, pages 133--142. IEEE, 2006.


\bibitem{Ziegler:LecturesPolytopes}
G.~M.~Ziegler, 
\newblock{Lectures on polytopes},
\newblock{Graduate Texts in Mathematics}, {152},
\newblock  {Springer-Verlag, 1995.}




\end{thebibliography}
